\documentclass{article}

\usepackage[a4paper,margin=2.5cm]{geometry}

\usepackage{color}
\usepackage[utf8]{inputenc} 
\usepackage[T1]{fontenc}    
\usepackage{hyperref}       
\usepackage{url}            
\usepackage{booktabs}       
\usepackage{amsfonts}       
\usepackage{nicefrac}       
\usepackage[expansion=false]{microtype}      
\usepackage{lipsum}
\usepackage{graphicx}
\usepackage{color}
\usepackage{amsthm,amsmath,amssymb}
\usepackage[english]{babel}
\title{On $\delta$-derivations of $\mathfrak{saff}_{n+1}$}

\author{
 Kazbek Nuratdinov\\
  National University of Uzbekistan\\
Tashkent, Uzbekistan\\
kazbeknur11@gmail.com\\
 \and
 Rustamjon Gaybullaev\\
 National University of Uzbekistan\\
Tashkent, Uzbekistan\\
r.gaybullaev@nuu.uz\\
 }

\date{}
\begin{document}
\maketitle
\begin{abstract}
We study $\delta$-derivations of the finite-dimensional perfect
Lie algebra
$\mathfrak{saff}_{n+1}
=\mathfrak{sl}_{n+1}(\mathbb C)\ltimes\mathbb C^{n+1}$,
where $\mathbb C^{n+1}$ is the natural
$\mathfrak{sl}_{n+1}(\mathbb C)$-module.
For $n>1$ and $\delta\ne1$, we prove that every nonzero
$\delta$-derivation occurs at $\delta=\frac12$
and is a scalar multiple of the identity.
For $n=1$, we explicitly describe the four-dimensional
space of $(-2)$-derivations.
\end{abstract}
\label{firstpage}
{\bf Keywords:} 
{Lie algebras, $\delta$-derivation.}

{\bf MSC 2020:} 17B40, 17B63

{\bf \large 1. Introduction}


Derivations are fundamental tools in the study of Lie algebras.
The notion of a $\delta$-derivation, introduced by Filippov
[2,3], generalizes the usual derivation identity by incorporating
a scalar parameter. A linear map $\varphi:L\to L$ is called
a $\delta$-derivation if
\[
\varphi([x,y])
=\delta\bigl([\varphi(x),y]+[x,\varphi(y)]\bigr)
\]
for all $x,y\in L$. The case $\delta=1$ gives ordinary
derivations, while every element of the centroid is
a $\frac12$-derivation. Thus, this notion provides a common
framework for studying derivations and centroidal maps.

Filippov established foundational results on
$\delta$-derivations of Lie algebras and prime Lie algebras
[2,3]. Zusmanovich [4] developed the theory further,
studying current Lie algebras, the modular Zassenhaus
algebra, and Lie superalgebras. In particular, he constructed
examples of $\frac12$-derivations that are zero divisors,
giving a negative answer to a question posed by Filippov.

Finite-dimensional perfect Lie algebras provide a natural
setting for extending these investigations beyond simple
and semisimple Lie algebras. A perfect Lie algebra satisfies
$[L,L]=L$, but it may have a nonzero solvable radical.
Consequently, the study of its derivations involves both
the structure of a Levi subalgebra and its action on
the radical. In their work on post-Lie algebra structures,
Burde and Dekimpe [1] also investigated derivations of
perfect Lie algebras, including semidirect products
with abelian radicals.

Related questions have been considered for other
non-semisimple Lie algebras. Yang and Tang [6] studied
derivations of the Schr\"odinger algebra and their
applications, while Tao [5] investigated simple modules
of the $n$-th Schr\"odinger algebra. These works illustrate
the interest in both derivations and representations
of Lie algebras with nontrivial semidirect-product structure.

In this paper, we consider the special affine Lie algebra
\[
\mathfrak{saff}_{n+1}
=\mathfrak{sl}_{n+1}(\mathbb C)\ltimes\mathbb C^{n+1},
\]
where $\mathbb C^{n+1}$ is the natural
$\mathfrak{sl}_{n+1}(\mathbb C)$-module.
This Lie algebra is perfect and has the abelian radical
$\mathbb C^{n+1}$. Its representation theory was studied
by Casati [7], who classified its finite-dimensional
perfect cyclic modules. Here we investigate
$\delta$-derivations of the same Lie algebra.

For $n>1$ and $\delta\ne1$, we give a complete description
of $\operatorname{Der}_{\delta}(\mathfrak{saff}_{n+1})$.
We prove that every $\frac12$-derivation is a scalar
multiple of the identity and that
$\operatorname{Der}_{\delta}(\mathfrak{saff}_{n+1})=0$
for $\delta\notin\{1,\frac12\}$.
The proof proceeds by establishing invariance of
the abelian radical, passing to the simple quotient,
and analyzing the remaining components using
the natural module.
Finally, we explicitly describe the four-dimensional
space of $(-2)$-derivations of $\mathfrak{saff}_2$,
which exhibits exceptional behaviour in the case $n=1$.

{\bf \large 2. Preliminaries}

For the convenience of the reader, we recall the relevant definitions and established results concerning $\delta$-derivations. Throughout this paper, $\mathbb{C}$ denotes the field of complex numbers.

\textbf{Definition 2.1.}
Let $(L,[\cdot ,\cdot])$ be an algebra with multiplication $[\cdot ,\cdot]$. Let $\varphi: L\longrightarrow L$ be a linear map and $\delta\in\mathbb{C}$. Then $\varphi$ is a $\delta$-derivation if it satisfies
\begin{equation}\label{12der}
\varphi\bigl([x,y]\bigr)=\delta\bigl([\varphi(x),y]+[x,\varphi(y)]\bigr).
\end{equation}

\textbf{Definition 2.2.}
The perfect Lie algebra $\mathfrak{saff}_{n+1}=\mathfrak{sl}_{n+1}(\mathbb{C})\ltimes \mathbb{C}^{n+1}$ is the Lie algebra with basis $\{h_{r},e_{ij},p_{k}\}$
whose non-trivial Lie brackets are:
\begin{align*}
[h_{r},e_{ij}]&=(\varepsilon_{ri}-\varepsilon_{rj}-\varepsilon_{r+1,i}+\varepsilon_{r+1,j})e_{ij},\\
[e_{ij},e_{kl}]&=\varepsilon_{jk}e_{il}-\varepsilon_{li}e_{kj},\ (j,i)\neq(k,l)\\
[h_{r},p_{k}]&=(\varepsilon_{rk}-\varepsilon_{r+1,k})p_{k},\\
[e_{ij},p_{k}]&=\varepsilon_{jk}p_{i},\\
[e_{ij},e_{ji}]&=\sum_{r=i}^{j-1}h_{r} \ (i<j),
\end{align*}
where $r\in \{1,\dots,n\},\ i,j,k,l\in\{1,\dots,n+1\}$ and $\varepsilon_{ij}$ is the Kronecker delta. All remaining brackets are zero or follow by skew-symmetry.

{\bf \large 3. Main results}

The following result gives a complete description of the space of $\delta$-derivations for $\mathfrak{saff}_{n+1}=\mathfrak{sl_{n+1}}(\mathbb{C})\ltimes \mathbb{C}^{n+1}$. In this work, we exclude the case $\delta=1$, corresponding to ordinary derivation. The case $\delta=1$ follows from Proposition~5.11 of
Burde and Dekimpe in [1], together with Schur's lemma applied to the standard $\mathfrak{sl}_{n+1}(\mathbb{C})$-module. Moreover, if $\delta=0$, then $\varphi([\mathfrak{saff}_{n+1},\mathfrak{saff}_{n+1}])=0$. Since $[\mathfrak{saff}_{n+1},\mathfrak{saff}_{n+1}]=\mathfrak{saff}_{n+1}$, it follows that $\varphi=0$. Henceforth, assume that $\delta\neq 0,1$. For convenience, set $\mathfrak{s}=\mathfrak{sl}_{n+1}(\mathbb{C})$
and $V=\mathbb{C}^{n+1}$. Then $\mathfrak{saff}_{n+1}=\mathfrak{s}\ltimes V$.
Before proving the main result, we first establish the following auxiliary lemma.

\textbf{Lemma 3.1.}
\emph{Every $\delta$-derivation $\varphi$ of $\mathfrak{saff}_{n+1}$ maps $V$ into itself, i.e. $\varphi(V)\subseteq V$.}

\medskip
\textbf{Proof.} First assume that $n>1$. Set

\begin{equation}
\varphi(p_{i})=\sum_{t=1}^{n} \alpha_{t}^{i}h_{t}+\sum_{p,q=1}^{n+1}\beta_{p,q}^{i}e_{p,q}+
\sum_{k=1}^{n+1}\gamma_{k}^{i}p_{k}, \ p\neq q,
\end{equation}
where $1\leq i\leq n+1$.
Applying (1) to $[h_{i},p_{i}]=p_{i}$ for $1\leq i<n+1$ and using (2), we obtain
\begin{equation}
\begin{aligned}
\sum_{t=1}^{n} \alpha_{t}^{i}h_{t}+\sum_{p,q=1}^{n+1}\beta_{p,q}^{i}e_{p,q}+
\sum_{k=1}^{n+1}\gamma_{k}^{i}p_{k}&=\delta\Bigl([\varphi(h_{i}),p_{i}] +\sum_{t=1}^{n} \alpha^{i}_{t} [h_{i},h_{t}]+\\
&+\sum_{p,q=1}^{n+1}\beta^{i}_{p,q}[h_{i},e_{p,q}]+\sum_{k=1}^{n+1}\gamma^{i}_{k}[h_{i},p_{k}]\Bigr).
\end{aligned}
\end{equation}
Since $V$ is an ideal, $[\varphi(h_{i}),p_{i}]$ and $[h_{i},p_{k}]$ lie in $V$, as $[h_{i},h_{t}]=0$ and contribute no $h_{i}$ component on the right-hand side, which implies $\alpha^{i}_{t}=0$ for $1\leq t<n+1$ in (3). Furthermore, applying (1) to $[h_{n},p_{n+1}]=-p_{n+1}, \ (i=n+1)$ we get an analogous result. For $1\leq i\leq n$, comparing the coefficients of $e_{p,q}$
in (3), we obtain $(1-\delta\tau)\beta^i_{p,q}=0$,
where $[h_i,e_{p,q}]=\tau e_{p,q}$. For $i=n+1$, applying (1) to
$[h_n,p_{n+1}]=-p_{n+1}$ and comparing the coefficients
of $e_{p,q}$ gives
$(-1-\delta\tau)\beta^{n+1}_{p,q}=0$. In both cases, $\tau\in\{0,\pm1,\pm2\}$.
Thus, a nonzero coefficient $\beta^i_{p,q}$ can occur only
if $\delta\tau=\pm1$, which implies $\delta=-1,-\frac{1}{2},\frac{1}{2},1$.
Since the case $\delta=1$ has already been excluded,
we conclude that $\beta^i_{p,q}=0$ whenever
$\delta\neq-1,-\frac{1}{2},\frac{1}{2}$. From (3) for $\delta=\frac12$ we arrive at
\begin{align}
\varphi(p_{i})&=\beta_{i,i+1}^{i}e_{i,i+1}+\sum_{k=1}^{n+1}\gamma^{i}_{k}p_{k}, \ 1\leq i\leq n,\\
\varphi(p_{n+1})&=\beta_{n+1,n}^{n+1}e_{n+1,n}+\sum_{k=1}^{n+1}\gamma^{n+1}_{k}p_{k}, \ i=n+1.
\end{align}
 Consider the relation $[h_{i-1},p_{i}]=-p_{i}$. Applying (1) to this relation and combining the result with (4), we obtain $\beta_{i,i+1}^{i}=0$ for $1<i<n+1$. If $i=1$ or $i=n+1$, it is sufficient to consider the relations $[h_{2},p_{1}]=0$ and $[h_{n-1},p_{n+1}]=0$, respectively. In the case $\delta=-\frac12$, it follows from (3) that
\begin{align}
\varphi(p_{i})&=\beta^{i}_{i+1,i}e_{i+1,i}+\sum^{n+1}_{k=1}\gamma^{i}_{k}p_{k}, \ 1\leq i\leq n,\\
\varphi(p_{n+1})&=\beta^{n+1}_{n,n+1}e_{n,n+1}+\sum^{n+1}_{k=1}\gamma^{n+1}_{k}p_{k}, \ i=n+1.
\end{align}
Likewise, applying (1) to the relation $[h_{i-1},p_{i}]=-p_{i}$ and together with (6) gives $\beta^{i}_{i+1,i}=0$ for $1<i\leq n$. Applying (1) to the relations $[h_{2},p_{1}]=0$ and $[h_{n-1},p_{n+1}]=0$ we obtain $\beta^{1}_{2,1}=0$ and $\beta^{n+1}_{n,n+1}=0$ for $i=1, n+1$, respectively. It remains to prove Lemma 3.1 for $\delta=-1$.  We first consider the case $1<i<n+1$, the boundary cases $i=1$ and $n+1$ will be shown separately. We rewrite (2) as follows:
\begin{align}
\varphi(p_{i})=\sum_{p,q=1}^{n+1}\beta_{p,q}^{i}e_{p,q}+\sum_{k=1}^{n+1}\gamma_{k}^{i}p_{k}, \ p\neq q.
\end{align}
Applying (1) to $[h_{i},p_{i}]=p_{i}$ together with (8) we obtain
\begin{equation}
\begin{aligned}
\sum_{p,q=1}^{n+1}\beta_{p,q}^{i}e_{p,q}+\sum_{k=1}^{n+1}\gamma_{k}^{i}p_{k}&=-\bigl([\varphi(h_{i}),p_{i}]+[h_{i},\varphi(p_{i})]\bigr)\\
&=-\bigl([\varphi(h_{i}),p_{i}]+\sum_{p,q=1}^{n+1}\beta_{p,q}^{i}[h_{i},e_{p,q}]+\sum_{k=1}^{n+1}\gamma_{k}^{i}[h_{i},p_{k}]\bigr).
\end{aligned}
\end{equation}
Comparing the corresponding coefficients at the basis elements $e_{p,q}$ in (9) yields $\beta^{i}_{p,q}=0$ with the exception of $i=q, i+1\neq p$ and $i+1=p, i\neq q$. That is,

\begin{align}
\varphi(p_{i})=\sum_{p=1}^{n+1}\beta_{p,i}^{i}e_{p,i}+\sum_{q=1}^{n+1}\beta_{i+1,q}^{i}e_{i+1,q}+\sum_{k=1}^{n+1}\gamma_{k}^{i}p_{k},
\ p,q \neq i,i+1.
\end{align}
Further, applying (1) to $[h_{i-1},p_{i}]=-p_{i}$ and using (8) gives
\begin{equation}
\begin{aligned}
-\sum_{p,q=1}^{n+1}\beta_{p,q}^{i}e_{p,q}-\sum_{k=1}^{n+1}\gamma_{k}^{i}p_{k}&=-\bigl([\varphi(h_{i-1}),p_{i}]+[h_{i-1},\varphi(p_{i})]\bigr)\\
&=-\bigl([\varphi(h_{i-1}),p_{i}]+\sum_{p,q=1}^{n+1}\beta_{p,q}^{i}[h_{i-1},e_{p,q}]+\sum_{k=1}^{n+1}\gamma_{k}^{i}[h_{i-1},p_{k}]\bigr).
\end{aligned}
\end{equation}
By the same argument from (11) we have the following
\begin{align}
\varphi(p_{i})=\sum_{q=1}^{n+1}\beta_{i-1,q}^{i}e_{i-1,q}+\sum_{p=1}^{n+1}\beta_{p,i}^{i}e_{p,i}+\sum_{k=1}^{n+1}\gamma_{k}^{i}p_{k},
\ p,q \neq i-1,i \ .
\end{align}
Comparing the corresponding coefficients at the basis elements in (10) and (12) we obtain $\beta^{i}_{i+1,q}=\beta^{i}_{i-1,q}=0$, more precisely
\begin{align}\label{eq:pi}
\varphi(p_{i})=\sum_{p=1}^{n+1}\beta_{p,i}^{i}e_{p,i}+\sum_{k=1}^{n+1}\gamma_{k}^{i}p_{k}, \ p\neq i-1,i,i+1.
\end{align}
Our final task is to show that $\beta_{p,i}^{i}=0$. Since for any $e_{pi}$ with $p\neq i-1,i,i+1$ there exists $h_{k} \ (k\neq i-1,i)$ such that $[h_{k},p_{i}]=0$ and $[h_{k},e_{pi}]\neq 0$, applying (1) to \([h_k,p_i]=0\) together with (13) implies $\beta_{p,i}^{i}=0$ immediately.\\
It only remains to show $\varphi(p_1)\in V, \ \varphi(p_{n+1})\in V.$ From (10), if $i=1$ we have
\begin{align}
\varphi(p_{1})=\sum_{p=3}^{n+1}\beta_{p,1}^{1}e_{p,1}+\sum_{q=3}^{n+1}\beta_{2,q}^{1}e_{2,q}+\sum_{k=1}^{n+1}\gamma_{k}^{1}p_{k}.
\end{align}
Applying (1) to $[h_{p},p_{1}]=0$ and $[h_{q},p_{1}]=0 \ (p,q\in \{3,\dots,n\})$ together with (14) we get $\beta_{p,1}^{1}=\beta_{2,q}^{1}=0$. If $p,q=n+1$, applying (1) to $[h_{n},p_{1}]=0$ we obtain $\beta_{n+1,1}^{1}=\beta_{2,n+1}^{1}=0$.\\
If $i=n+1$, from (12) we have
\begin{align}
\varphi(p_{n+1})=\sum_{q=1}^{n-1}\beta_{n,q}^{n+1}e_{n,q}+\sum_{p=1}^{n-1}\beta_{p,n+1}^{n+1}e_{p,n+1}+\sum_{k=1}^{n+1}\gamma_{k}^{n+1}p_{k}.
\end{align}
Applying (1) to $[h_{q},p_{n+1}]=0$ and $[h_{p},p_{n+1}]=0$ together with (15) we obtain $\beta_{n,q}^{n+1}=\beta_{p,n+1}^{n+1}=0$ as desired.\\
Let $n=1$, then $\mathfrak{saff}_{2}=\operatorname{span}\{h_{1},e_{12},e_{21},p_{1},p_{2}\}$. Write
\begin{align*}
\varphi(p_{1})&=\alpha_{1}h_{1}+\beta_{1}e_{12}+\gamma_{1}e_{21}+\xi_{1}p_{1}+\mu_{1}p_{2},\\
\varphi(p_{2})&=\alpha_{2}h_{1}+\beta_{2}e_{12}+\gamma_{2}e_{21}+\xi_{2}p_{1}+\mu_{2}p_{2},
\end{align*}
where $h_{1},e_{12},e_{21}\in\mathfrak{s}, p_{1},p_{2}\in V$. Applying (1) to $[e_{12},p_{1}]=0,[e_{21},p_{2}]=0$ and $[p_{1},p_{2}]=0$, we establish $\varphi(V)\subseteq V$, completing the proof.

\textbf{Remark 3.2.} Note that if $n=2$, then we immediately get $\beta_{p,i}^i=0$ from (13).

In what follows, we assume that $n>1$. Let $\mathfrak{saff}_{n+1}=\mathfrak{s}\ltimes V$, where $V$ is regarded as the natural $\mathfrak{s}$-module. Recall that $V$ is an abelian ideal of $\mathfrak{saff}_{n+1}$. We consider the quotient Lie algebra $\mathfrak{saff}_{n+1}/V$. Since, as a vector space, $\mathfrak{saff}_{n+1}=\mathfrak{s}\oplus V$, every coset in $\mathfrak{saff}_{n+1}/V$ can be represented by an element of $\mathfrak{s}$. More precisely, for
$x\in\mathfrak{s}$ and $v\in V$, one has $(x+v)+V=x+V$.
The Lie bracket in the quotient is given by $[x'+V,y'+V]=[x',y']+V$. Let $x'=x+u,\ y'=y+v$, where $x,y\in\mathfrak{s}$ and $u,v\in V$. By the definition of the semidirect product, $[x',y']=[x,y]+x\cdot v-y\cdot u$. Since $x\cdot v-y\cdot u\in V$, it follows that $[x',y']+V=[x,y]+V$. Therefore, $[x'+V,y'+V]=[x,y]+V$. Consequently, the map
\[
\Phi:\mathfrak{saff}_{n+1}/V\longrightarrow \mathfrak{s}, \ \Phi(x+V)=x,
\]
is a well-defined Lie algebra isomorphism. Hence, $\mathfrak{saff}_{n+1}/V\cong\mathfrak{s}$.
Now, let $\varphi:\mathfrak{saff}_{n+1}\to\mathfrak{saff}_{n+1}$ be a $\delta$-derivation satisfying $\varphi (V)\subseteq V$. Then $\varphi$ induces a linear map $\tilde{\varphi}:\mathfrak{saff}_{n+1}/V\to\mathfrak{saff}_{n+1}/V$ defined by
\[
\tilde{\varphi}(x'+V)=\varphi(x')+V.
\]
The condition $\varphi (V)\subseteq V$ guarantees that $\tilde{\varphi}$ is well defined. Moreover, $\tilde{\varphi}$ is again a $\delta$-derivation. Indeed, for arbitrary $x',y'\in\mathfrak{saff}_{n+1}$,
\begin{align*}
\tilde{\varphi}\bigl([x'+V,y'+V]\bigr)&=\varphi\bigl([x',y']\bigr)+V\\
&=\delta\bigl([\varphi(x'),y']+[x',\varphi(y')]\bigr)+V\\
&=\delta\bigl([\tilde{\varphi}(x'+V),y'+V]+[x'+V,\tilde{\varphi}(y'+V)]\bigr).
\end{align*}
Thus, $\tilde{\varphi}\in\operatorname{Der}_{\delta}(\mathfrak{saff}_{n+1}/V)$.
For $\delta\notin\{-1,0,\tfrac12,1\}$, we have
$\widetilde{\varphi}=0$ by [3, Theorem 2].
If $\delta=-1$, then $\widetilde{\varphi}=0$
by [3, Corollary 1], since $n>1$ and hence
$\mathfrak{sl}_{n+1}(\mathbb C)$ is not isomorphic
to $\mathfrak{sl}_2(\mathbb C)$.
Thus, under the standing assumption $\delta\ne0,1$,
we obtain $\widetilde{\varphi}=0$ for
$\delta\ne\tfrac12$.
Hence, $\varphi(\mathfrak{saff}_{n+1})\subseteq V$. Now, we show that $\varphi(\mathfrak{s})=0$ for
$\delta\neq \frac12,1$ (computation (16)--(23) below).
If $\delta=\frac12$, then $\tilde{\varphi}=\lambda\operatorname{Id}$
for some $\lambda\in\mathbb{C}$, since every
$\frac12$-derivation of the simple Lie algebra
$\mathfrak{s}$ is a scalar multiple
of the identity map (see [2, Corollary~3]).
Then $\psi=\varphi-\lambda\operatorname{Id}$ is a
$\frac12$-derivation of $\mathfrak{saff}_{n+1}$ with
$\psi(V)\subseteq V$
and $\tilde{\psi}=0$, i.e.
$\psi(\mathfrak{s})\subseteq V$.
Since (16)--(23) use only $\delta\neq 0,1$ and $n>1$,
the same computation applied to $\psi$ gives
$\psi(\mathfrak{s})=0$, i.e.
$\varphi(x)=\lambda x$ for all $x\in\mathfrak{s}$. \\
For the following computation, assume that
$\varphi(\mathfrak s)\subseteq V$ and $\delta\ne0,1$. Let   $H\subset\mathfrak{s}$ be the Cartan subalgebra. For $h,k\in H,$ set
\begin{equation}
\begin{aligned}
\varphi(h)&=\sum_{r=1}^{n+1}\alpha_{r}(h)p_{r}, \ h=diag(h_{1},\dots,h_{r},\dots,h_{n+1}),\\
\varphi(k)&=\sum_{r=1}^{n+1}\alpha_{r}(k)p_{r}, \ k=diag(k_{1},\dots,k_{r},\dots,k_{n+1}),
\end{aligned}
\end{equation}
where $tr(h)=0$ for all $h\in H$.\\
Applying (1) to $[h,k]=0$ for all $h,k\in H$ together with (16) we obtain
\[
[h,\sum_{r=1}^{n+1}\alpha_{r}(k)p_{r}]=[k,\sum_{r=1}^{n+1}\alpha_{r}(h)p_{r}],
\]
or equivalently
\[
\sum_{r=1}^{n+1}h_{r}\alpha_{r}(k)p_{r}=\sum_{r=1}^{n+1}k_{r}\alpha_{r}(h)p_{r},
\]
by comparing both sides we get $h_{r}\alpha_{r}(k)=k_{r}\alpha_{r}(h)$. Here we can choose $k\in H$ with $k_{r}\neq 0$ for $1\leq r\leq n+1$, keeping $tr(k)=0$. Consequently, there exist scalars $a_{r}\in\mathbb{C}$ such that $\alpha_{r}(h)=a_{r}h_{r}$. Hence,
\begin{align}
\varphi(h)=hp_{0}, \text{ where } p_{0}=\sum_{r=1}^{n+1}a_{r}p_{r}.
\end{align}
Let
\begin{align}
\varphi(e_{ij})=\sum_{r=1}^{n+1}b_{r}p_{r}, \ i\neq j.
\end{align}
Applying (1) to $[h,e_{ij}]=(h_{i}-h_{j})e_{ij}$ together with (17) and (18) we obtain
\begin{equation}
\begin{aligned}
(h_{i}-h_{j})\sum_{r=1}^{n+1}b_{r}p_{r}&=\delta\bigl([\varphi(h),e_{ij}]+[h,\varphi(e_{ij})]\bigr)\\
&=\delta\bigl([h\sum_{r=1}^{n+1}a_{r}p_{r},e_{ij}]+[h,\sum_{r=1}^{n+1}b_{r}p_{r}]\bigr)\\
&=\delta\bigl(\sum_{r=1}^{n+1}h_{r}b_{r}p_{r}-h_{j}a_{j}p_{i}\bigr),
\end{aligned}
\end{equation}
it follows from (19)
\begin{align}
(h_{i}-h_{j})b_{r}&=\delta h_{r}b_{r},\ r\neq i\\
(h_{i}-h_{j})b_{i}&=\delta(h_{i}b_{i}-h_{j}a_{j}), \ r=i.
\end{align}
It follows from (20)
\begin{align}
(h_{i}-h_{j}-\delta h_{r})b_{r}=0.
\end{align}
If $r\neq j$, choosing $h\in H$ with $h_{i}=1, h_{j}=-1$ and all remaining diagonal entries equal to zero, we get $b_{r}=0 \ (r\neq i,j)$. Further, if $r=j$, choosing $h_{i}=1, h_{j}=0, h_{s}=-1 \ (s\neq i,j)$ and $h_{l}=0$ everywhere else, then (22) gives $b_{j}=0$. Hence, $\varphi(e_{ij})=b_{i}p_{i}$. \\
On the other hand, rewrite (21) as follows
\begin{align}
(1-\delta)b_{i}h_{i}+(\delta a_{j}-b_{i})h_{j}=0.
\end{align}
Choosing $h_{i}=1, h_{j}=0, h_{s}=-1$, and then choosing $h_{i}=0, h_{j}=1, h_{s}=-1$ from (23) we get
\begin{align*}
(1-\delta)b_{i}=0, \ b_{i}=\delta a_{j},
\end{align*}
as $\delta\neq 0, 1$, whence $a_{j}=b_{i}=0$. For each $j$, choosing such an $i\neq j$, we obtain $a_{j}=0$ for $1\leq j\leq n+1$. Thus, $\varphi(\mathfrak{s})\equiv0$. \\
By Lemma 3.1, $\varphi(p_{j})\in V$. Write
\begin{equation}
\varphi(p_{j})=\sum_{t=1}^{n+1}\mu_{t}^{j}p_{t}, \ 1\leq j\leq n+1.
\end{equation}
Applying (1) to $[h_{j},p_{j}]=p_{j} \ (1\leq j\leq n)$ and using (24), we obtain
\begin{equation}
\begin{aligned}
\sum_{t=1}^{n+1}\mu_{t}^{j}p_{t}&=\delta\bigl([\varphi(h_{j}),p_{j}]+[h_{j},\varphi(p_{j})]\bigr)\\
&=\delta\bigl([\varphi(h_{j}),p_{j}]+\mu_{j}^{j}p_{j}-\mu_{j+1}^{j}p_{j+1}\bigr).
\end{aligned}
\end{equation}
We consider the cases $\delta=\frac{1}{2}$ and $\delta\neq\frac{1}{2}$ separately. If $\delta=\frac{1}{2}$, then $\varphi(h_{j})=\lambda h_{j}$ and comparing the corresponding coefficients at basis elements in (25), we obtain $\mu^{j}_{j}=\lambda$ and $\mu^{j}_{t}=0$ for $t\neq j$. Thus, $\varphi(p_{j})=\lambda p_{j}$ for $1\leq j\leq n$. For $j=n+1$, applying (1) to $[h_{n},p_{n+1}]=-p_{n+1}$ gives $\varphi(p_{n+1})=\lambda p_{n+1}$.\\
Now suppose that $\delta\neq \frac{1}{2}$. Applying (1) to $[h_{j},p_{j}]=p_{j}$ for $1\leq j\leq n$ we obtain
\begin{align*}
\sum_{t=1}^{n+1}\mu_t^j p_t=\delta\left(\mu_j^j p_j-\mu_{j+1}^j p_{j+1}\right).
\end{align*}
Consequently,
\begin{align*}
(1-\delta)\mu_j^j=0, \ (1+\delta)\mu_{j+1}^j=0, \ \mu_t^j=0 \quad t\neq j,j+1.
\end{align*}
Similarly, the relation $[h_n,p_{n+1}]=-p_{n+1}$
implies
\begin{align*}
(1-\delta)\mu_{n+1}^{n+1}=0, \ (1+\delta)\mu_n^{n+1}=0, \ \mu_t^{n+1}=0 \quad t\neq n,n+1.
\end{align*}
Thus, if $\delta\neq-1$, then
$\varphi(p_j)=0$ for every $j$.\\
It remains to consider $\delta=-1$. The preceding
equalities reduce to
\begin{align*}
\varphi(p_j)=\mu_{j+1}^j p_{j+1} \ (1\leq j\leq n), \ \varphi(p_{n+1})=\mu_n^{n+1}p_n.
\end{align*}
For $2\leq j\leq n$, applying (1) with $\delta=-1$ to $[h_{j-1},p_j]=-p_j$ and using
$[h_{j-1},p_{j+1}]=0$ gives $\mu_{j+1}^j=0$.
For the remaining indices, the relations
\begin{align*}
[h_2,p_1]=0, \ [h_{n-1},p_{n+1}]=0
\end{align*}
give $\mu_2^1=0$ and $\mu_n^{n+1}=0$, respectively.

Together with Lemma 3.1 and the calculations above, we obtain the following theorem, which is the main result of this paper.\\

\textbf{Theorem 3.3.}
\emph{Let $n>1$ and $\delta\in\mathbb C\setminus\{1\}$.
Then
\[
\operatorname{Der}_{\delta}(\mathfrak{saff}_{n+1})=
\begin{cases}
\mathbb{C}\cdot\operatorname{Id},
& \delta=\frac{1}{2},\\[2mm]
0,
& \delta\notin\left\{\frac{1}{2},1\right\}.
\end{cases}
\]}

\textbf{Remark 3.4.}
\emph{The case $n=1$ is exceptional. In this case, $\mathfrak{saff}_{2}$ admits nontrivial $(-2)$-derivations. More precisely, let $\mathfrak{saff}_{2}=\operatorname{span}\{h_{1},e_{12},e_{21},p_{1},p_{2}\}$, then every $(-2)$-derivation is given by four arbitrary parameters $\alpha,\beta,\gamma,\xi\in \mathbb{C}$ through
\begin{align*}
\varphi(h_{1})&=-2\xi p_{1}+2\alpha p_{2}, \varphi(e_{12})=\alpha p_{1}+\beta p_{2},\varphi(e_{21})=\gamma p_{1}+\xi p_{2},\\
\varphi(p_{1})&=0, \ \varphi(p_{2})=0.
\end{align*}}



\medskip



\medskip



\medskip



\medskip


\medskip





\medskip


\newpage



\medskip





\bigskip

\textbf{Received: **/**/202*}\\

\textbf{Accepted: **/**/202*}\\

\bigskip
\textcolor{red}{\Large Cite this article}

K.~Nuratdinov, R.~Gaybullaev: {On $\delta$-derivations of $\mathfrak{saff}_{n+1}$.}  \textit{Bull.~Inst.~Math.},~ 2023, Vol.6, No 1, pp. \pageref{firstpage}-\pageref{lastpage}

\label{lastpage}

\end{document}